\documentclass[a4paper,12pt]{article}
\usepackage{amsmath,amsfonts}
\usepackage[english]{babel}
\usepackage{a4wide}

\usepackage{amssymb,amsthm}

\usepackage{color}

\usepackage[cp1251]{inputenc}

\theoremstyle{plain}

\newtheorem{thm}{Theorem}
\newtheorem{lemma}{Lemma}
\newtheorem{cor}{Corollary}
\newtheorem{remark}{Remark}
\theoremstyle{definition}

\title{On asymptotic behavior of almost surely extreme values of independent random variables }

\author{Kateryna Akbash$^a$, Ivan Matsak$^b$}

\date{}

\begin{document}
\maketitle

\begin{center}
$^a$ Volodymyr Vynnychenko Central Ukrainian State University,
Shevchenko street 1, Kropyvnytskyi 25006, Ukraine \\

$^b$ Taras Shevchenko National University of Kyiv,
$2/6 $, Academician Glushkov avenue, Kyiv 03127, Ukraine\\
\medskip
kateryna.akbash@gmail.com, ORCID 0000-0003-3676-4574 (K.Akbash), ivanmatsak@univ.kiev.ua, ORCID 0000-0003-3816-7650 (I.Matsak)

\end{center} \bigskip

\hspace{5,0cm}

\noindent{\it Keywords:}  extreme values,   discrete random
variables,  almost sure limit theorems

\medskip

\noindent{\it MSC:} Primary 60G70,  60F15

\medskip

\bigskip

\begin{quote}

\begin{center}

                   \vspace{0.5cm}

{Abstract}\\
The article studies the almost surely asymptotics of extreme values $\bar{\xi}_n = \max_{1\leq i \leq n} \xi_i$, where $ \xi , \xi_1 , \xi_2 , \ldots$ are discrete identically distributed random variables. One of the main results on this topic is related to the law of the iterated logarithm for the lim sup (LIL) and a law of the triple logarithm for the lim inf (LTL). But, taking into account the specifics of the discrete case, necessary and sufficient conditions are established for
$\mathbf{P}(\bar{\xi}_n = a_n +l \quad\mbox{infinitely often} ) =1$, where $a_n$ is some given increasing sequence of integers and $l$ is a fixed integer. Note that in the case of discrete random variables whose distribution tails are close to the tails of the Poisson distribution or fall off even faster, Theorems 1-3 of this article are significantly more informative than LIL and LTL. For the geometric distribution (or discrete random variables whose tails fall off even more slowly), the results are given in Theorem 3 and Theorem B, which are an important complement to LIL and LTL.

\end{center}

\end{quote}

\bigskip

\maketitle

\section{  Introduction. The main results.}

 Let's consider the sequence  $ \xi , \xi_1 , \xi_2 , \ldots$    of independent identically distributed random variables (i.i.d.r.v.)  with a distribution function (d.f.)  $F(x)=
\mathbf{P}(\xi < x ) $. And let
\begin{equation} \label{f1}
\bar{\xi}_n = \max_{1\leq i \leq n} \xi_i .
\end{equation}
The importance of the limit theorems for r.v. $\bar{\xi}_n $  is
clear. Therefore, the asymptotic of $\bar{\xi}_n $  when $n\rightarrow\infty$  has been studied extensively (see.\cite{ha_06}, \cite{gal}, \cite{gl_kou}, \cite{gne},  \cite{llr}). However, the weak convergence of normalized random variable r.v. $\bar{\xi}_n $  has been studied. At the same time, the problems on asymptotic behavior of almost surely (a.s.) extreme values of
i.i.d.r.v. haven't been properly studied.

Here we will be interested in just such problems. One of the important areas of research on this topic is related to the law of the iterated logarithm for the lim sup (LIL) and a law of the triple logarithm for the lim inf (LTL).

The law of the iterated logarithm for extreme values $\bar{\xi}_n $ of normally distributed r.v. $\xi_i$ was first considered in the article by J.Pickands \cite{pic}.


The next step was taken in the article by J. Kiefer \cite{kie72} (see also the article by J.~Fill and D.~Naiman \cite{fil20}).

For the case when the r.v. $\xi$ has an exponential distribution with the parameter $\lambda =1$, the following equalities were obtained in \cite{kie72}:
\begin{equation} \label{f1*}
\mathbf{P}\left(\bar{\xi}_n \geq L(n)+cL_2(n) \quad \mbox{i.o.}\right)=
\left\{
\begin{array}{l}
1, \quad \mbox{if} \quad c \leq 1, \\
0, \quad \mbox{if} \quad c > 1,
\end{array}
\right.
\end{equation}

\begin{equation} \label{f2*}
\mathbf{P}\left(\bar{\xi}_n \leq L(n)-L_3(n) + c \quad \mbox{i.o.}\right)=
\left\{
\begin{array}{l}
1, \quad \mbox{if} \quad c \geq 0, \\
0, \quad \mbox{if} \quad c < 1,
\end{array}
\right.
\end{equation}
where $L(n)=L_1(n)=\log n$, $L_k(n)=\log L_{k-1}(n)$ and $"\mbox{i.o.}"$ means infinitely often.

From equalities (\ref{f1*}) and (\ref{f2*}), LIL and LTL immediately
follow:
$$\limsup_{n\rightarrow\infty }{\frac{\bar{\xi}_n - L(n)}{L_2(n)}}=1 \quad \mbox{a.s.} \quad \mbox{and} \quad \liminf_{n\rightarrow\infty }{\frac{\bar{\xi}_n - L(n)}{L_3(n)}}=-1 \quad \mbox{a.s.}$$
(see also article by P. Glasserman et al. \cite{gl_kou}).

It should also be noted that articles \cite{kie72} and \cite{fil20}
investigate the general case of r.v. $(\xi_i)$ that take values in
$\mathbb{R}^d$, where $d\geq1$.

For further results concerning LIL and LTL in the case of more general random variables and some regenerative random processes see \cite{adm21}, \cite{adm24}, \cite{am2}. These publications provide some overview of the results on this topic, but they are far from complete. In more detail, the early results on the asymptotic behavior of extreme values of i.i.d.r.v., which were published before 1978, are considered in Chapter 4 of Book \cite{gal}.

In this article we will continue the study of the case of discrete r.v., which was started in articles  \cite{mi16}, \cite{mi_tims19}. Here we will obtain the necessary and sufficient conditions for almost surely:
$$(\bar{\xi}_n =a_n+l \quad \mbox{i.o.}),$$
where $a_n$ is an increasing sequence of integers that is defined in equality (\ref{f2}), $l$ is a fixed number.

It is well known that the asymptotics of extreme values in the continuous and discrete cases can differ significantly. This is confirmed by the following theorems 1-3. In our opinion, the importance of the obtained results is as follows. In the case of discrete random variables, the tails of the distribution are close to the tails of the Poisson distribution or fall off even faster, and theorems 1-3 and their conclusions are significantly more informative than LIL and LTL.

If we consider a geometric distribution (or discrete r.v., whose tails decay even more slowly), the situation changes slightly. But here too the results given in Theorem 3 and Theorem B are an important addition to LIL and LTL.


Thus, we will study the discreet case, more precisely, we assume
that r.v.  ${\xi} $  has a distribution
  \quad $(i , p_i ), i\geq 0, $ :
\[
\mathbf{P}(\xi =i ) = p_i > 0, \quad \sum_{i=0}^{\infty} p_i =1 .
\]
For such random variable $\xi$ we introduce the notation:
$$R({n})= -\ln\mathbf{P}(\xi  \geq n )= -\ln\left(\sum_{i\geq n} p_i \right),$$
$$ r(n) = R(n) - R(n-1) ,$$
\begin{equation} \label{f2}
 a_{n}= \max\left( k\geq 0: \sum_{i\geq k} p_i \geq \frac{1}{n} \right).
\end{equation}

Let us start with the result of the article  \cite{mi16} that
describes the asymptotic a.s. of the extreme values of discrete
 r.v., which tails of distributions decrease rather
quickly. The following theorem is given in \cite{mi16} in a slightly
different but equivalent formulation.

\textbf{Theorem  A.} \emph{Let } $\xi$ \emph{be discrete r.v.  with distribution} $(i , p_i )$, $i\geq 0$, $\bar{\xi}_n$, $r(n)$ and $a_n $ \emph{ are defined in
equalities} (\ref{f1}) and (\ref{f2}) \emph{respectively},
$r(n) $  \emph{is monotonic function and }
\begin{equation} \label{f3}
r(n)\rightarrow\infty , \quad  n \rightarrow\infty .
\end{equation}

 \emph{Then }

(i) \emph{equality}
\begin{equation} \label{f4}
\mathbf{P}\left(\limsup_{n \rightarrow\infty} (\bar{\xi}_n -  a_n ) = 1 \right)
= 1
\end{equation}
\emph{holds if and only if}
\begin{equation*}
\sum_{k\geq 1}  \exp(- r(k))  < \infty;
\end{equation*}

(ii) \emph{ to satisfy the equality }
\begin{equation} \label{f5}
\mathbf{P}\left(\liminf_{n \rightarrow\infty} (\bar{\xi}_n -  a_n ) = -1
\right) = 1
\end{equation}
 \emph{it is necessary and
sufficient that}
\begin{equation*}
\sum_{k\geq 1}  \exp(- e^{r(k)})  < \infty .
\end{equation*}

\begin{remark}\label{r1}
a) In fact, Theorem A is a certain discrete analogue of the
well-known result of O.Barndorff-Nielsen \cite{bn} on the asymptotic
stability of  $\bar{\xi}_n $  for continuous r.v.

b)  It should be noted that in the case of fast tails, Theorem A
gives an almost complete description of the asymptotic behavior of
extreme values a.s. Thus, for example, under the condition
\[
r(k)= C  L(k) + o(L_2 (k)), \quad C>1,
\]
 equalities (\ref{f4}), (\ref{f5})   of Theorem A are holds   and
\[
\mathbf{P}(\bar{\xi}_n = a_n  \quad \mbox{i.o.} ) = 1,
\]
(see. \cite{adm21}).

\end{remark}

 Let us formulate the main results of the article.

\begin{thm}\label{t1}
 Let  $\xi$  be a discrete random variable with the
distribution  $(i , p_i )$, $i\geq 0$, $\bar{\xi}_n , r(n)$
and  $a_n $ are defined in equalities  (\ref{f1}) and
(\ref{f2}) respectively, $l $ is a fixed integer, $l \geq 0$. And let  $r(n) $ be a monotonic function for which condition (\ref{f3}) is satisfied.  Then

  (i) \quad for $ l\in \{0, 1\}$
\[
\mathbf{P}(\bar{\xi}_n = a_n +l  \quad \mbox{i.o.} )=1,
\]

(ii) \quad for  $ l>1$   probability
\begin{equation}\label{f6}
\mathbf{P}(\bar{\xi}_n = a_n +l  \quad \mbox{i.o.} )
\end{equation}
 is equal to zero or one according to whether the following series converges or diverges
\begin{equation}\label{f7}
\sum_{k\geq 1}  \exp(- j r(k))
\end{equation}
for $j=l-1$.

\end{thm}

\begin{cor}\label{c1}
 Let  $\xi$ { be a discrete r.v. with the distribution} $(i , p_i )$, $i\geq 0$, $\bar{\xi}_n , r(n)$ and $a_n $ are defined in equalities  (\ref{f1}) and (\ref{f2})  respectively, $m $ is a fixed integer, $m \geq1$. And let $r(n) $ be monotonic function for which the condition  (\ref{f3})  is hold.  Then if series  (\ref{f7}) converges for  $j=m$  and diverges for $j=m-1$, then

(i)
\begin{equation}\label{f8}
\mathbf{P}\left(\limsup_{n \rightarrow\infty} (\bar{\xi}_n -  a_n ) = m \right)= 1,
\end{equation}
and the equality (\ref{f5}) of Theorem A is satisfied;

(ii) $\forall \, l=-1, 0, 1,\ldots, m $
\begin{equation}\label{f9}
\mathbf{P}\left(\bar{\xi}_n = a_n +l  \quad \mbox{i.o.} \right) = 1.
\end{equation}

\end{cor}

\begin{thm}\label{t2}
  Let  $\xi$ {be discrete r.v.  with distribution} $(i , p_i )$, $i\geq 0$, $\bar{\xi}_n, r(n)$ and $a_n $ defined in equalities  (\ref{f1}) and (\ref{f2}) respectively, \quad  $l $ - fixed integer, $l \geq 1$. And let $r(n)$ be a monotonic function for which the condition  (\ref{f3})
is hold.  Then

  (i) for  $l=1 $
\begin{equation}\label{f10}
\mathbf{P}(\bar{\xi}_n \leq a_n -l  \quad \mbox{i.o.} ) =1;
\end{equation}

  (ii ) if $l>1$ , then the probability
\begin{equation}\label{f111}
\mathbf{P}(\bar{\xi}_n \leq a_n -l  \quad \mbox{i.o.} )
\end{equation}
  is equal to zero or one in accordance with whether the follow series converges or diverges
\begin{equation}\label{f11}
\sum_{k\geq 1}  \exp(- \exp(j r(k)))
\end{equation}
for  $j=l-1$.

\end{thm}

From Theorem  \ref{t2} it is clear that when the series  (\ref{f11})
 converges at $j=l-1$, then
\[
\mathbf{P}(\bar{\xi}_n = a_n -l  \quad \mbox{i.o.} )=0 .
\]
 But it is not clear under
what conditions
 \begin{equation}\label{f100}
\mathbf{P}(\bar{\xi}_n = a_n -l  \quad \mbox{i.o.} ) =1.
\end{equation}
The following result basically answers this question.

\begin{thm}\label{t3}
If, under the conditions and notation of Theorem  \ref{t2}, the function $r(n) $ satisfies the condition
\begin{equation}\label{f101}
\lim_{n\rightarrow\infty} \frac{r(n)}{n} =0
\end{equation}
 and diverges series (\ref{f11}) for $j=l-1$, then  equality  (\ref{f100}) is hold.

\end{thm}

\begin{cor}\label{c2}
 Let  $\xi$ be a discrete r.v.  with the distribution $(i , p_i )$, $i\geq 0$, $\bar{\xi}_n$, $r(n)$
and $a_n $ are defined in the equalities (\ref{f1}) and (\ref{f2}) respectively,   $m$ is a fixed integer, $m \geq 1$. And let $r(n)$  be a monotonic function for which the condition  (\ref{f3}) holds.

  Then if series
\begin{equation} \label{f12}
\sum_{k\geq 1}  \exp(- e^{jr(k)})
\end{equation}
 converges for $j=m $,  and diverges for $j=m-1 $, then

(i)
\begin{equation} \label{f13}
\mathbf{P}\left(\liminf_{n \rightarrow\infty} (\bar{\xi}_n -  a_n ) =-m\right) = 1;
\end{equation}

(ii) for $\forall \, l \geq 0 $ the equality  (\ref{f9}) holds, and
if the function $r(n) $ satisfies the condition  (\ref{f101}) then
$\forall \, l \in (m, m-1, \ldots, 1) $  the equality  (\ref{f100}) holds.

\end{cor}

In conclusion, for completeness, we present the result of the
article \cite{adm24}, which investigates the case when the tails
of the distribution of r.v. are either close to the geometric
distribution or decrease more slowly.

\textbf{Theorem  B.}  \emph{Let}  $\xi , \xi_1 , \xi_2 , \ldots$ \emph{be a sequence of discrete  i.i.d.r.v.,  with the
distribution} $(i , p_i )$, $i\geq 0$, $\bar{\xi}_n$, $r(n)$
\emph{and} $a_n $ \emph{are defined in the equalities} (\ref{f1}) and (\ref{f2}) \emph{ respectively}. \emph{If one of the following conditions holds:}

(i) \quad $r(n)$ \emph{is an increasing function and for} $n
\rightarrow\infty$
\begin{equation} \label{f120}
r(n)= \chi(n) L_2 (n) \rightarrow \infty, \quad \chi(n)\rightarrow 0
\end{equation}

\emph{or}

(ii)\quad $r(n)$ \emph{satisfies the condition}

 \[
\exists C_0 < \infty, \quad \forall n\geq 1  \quad r(n) \leq  C_0,
 \]

 \emph{then for any integer} $m$
\begin{equation} \label{f130}
\mathbf{P}(\bar{\xi}_n = a_n +m  \quad \mbox{i.o.} )  = 1.
\end{equation}

 In the following sections  $2 - 4$ we prove the above Theorems \ref{t1}-\ref{t3}.
Finally, we conclude with some examples of applications of the above
results to the study of the asymptotic of extreme values of random
variables .

\section{ Proof of Theorem  \ref{t1}.}

First, let's give some auxiliary Lemmas.
\begin{lemma}\label{l2.1}
Let  $A_1 , A_2 ,  \ldots \quad  $ is an infinite sequence of random events. Then

$(i)$ \quad {if }
$$\sum_{n=1}^{\infty} \mathbf{P}(A_n)<\infty,$$
then
$$\mathbf{P}\left(\bigcap_{n=1}^{\infty} \bigcup_{k=n}^{\infty} A_k\right) =0; $$

$(ii)$ \quad {if }
\begin{equation} \label{f2.1}
\sum_{n=1}^{\infty} \mathbf{P}(A_n)=\infty,
\end{equation}
and for some $K>0$
\begin{equation} \label{f2.2}
\liminf_{n\rightarrow\infty} \frac{\sum_{j=1}^{n}
\sum_{k=1}^{n}\mathbf{P}(A_j \bigcap
A_k)}{(\sum_{j=1}^{n}\mathbf{P}(A_j))^2} \leq K ,
\end{equation}
then
\begin{equation} \label{f2.3}
  \quad  \mathbf{P}\left(\bigcap_{n=1}^{\infty} \bigcup_{k=n}^{\infty} A_k\right)\geq \frac{1}{K}.
\end{equation}

\end{lemma}

Lemma \ref{l2.1} is some generalization of a well-known Borel-Cantelli Lemma (see \cite{spi},  Chapter 6, \S 26).

\begin{remark}\label{r2}
Let  $(\eta_n )$, $n\geq 1$, be a sequence of i.r.v., $\mathfrak{B}_n = \mathfrak{B}(\eta_n , \eta_{n+1}, \ldots )$ is the $\sigma$-algebra generated by r.v. $\eta_n, \eta_{n+1}, \ldots$, $\mathfrak{B} = \bigcap_{n=1}^{\infty} \mathfrak{B}_n $ is residual   $\sigma$-algebra. Suppose that the random events  $ A_n $ defined in Lemma \ref{l2.1}  depend on r.v. $(\eta_n )$ such that
 $\bigcap_{n=1}^{\infty} \bigcup_{k=n}^{\infty} A_k \in
\mathfrak{B}$ and the conditions (\ref{f2.1}), (\ref{f2.2}) are
satisfied. Then, according to the Hewitt-Savage zero-one law in conditions of Lemma
 \ref{l2.1} the value  $\mathbf{P}(\bigcap_{n=1}^{\infty}\bigcup_{k=n}^{\infty} A_k)$ can take only two values: $0$ or $1$. Consequently, given the estimate  (\ref{f2.3}), we get

$$\mathbf{P}\left(\bigcap_{n=1}^{\infty} \bigcup_{k=n}^{\infty} A_k\right)  =1. $$

\end{remark}

Next, for the sequence of discrete i.i.d.r.v. $\xi , \xi_1 , \xi_2 , \ldots,$  with distribution  $(i , p_i )$, $i\geq 0,$  we construct random events $A_n$  as follows:
\begin{equation} \label{f2.4}
 \quad A_n =\{ \xi_n =\bar{\xi}_n = a_n +l \} ,
\end{equation}
where $a_n$ is defined by formula  (\ref{f2}),  $l$ is some integer.

\begin{lemma}\label{l2.2}
 Let random events  $A_n$ be given by the equality (\ref{f2.4}), $l \geq 0$ is a fixed integer. And let $r(n) $ be a monotonic function defined in  (\ref{f2}),  for which the condition
 (\ref{f3}) is satisfied. Then

(i) \quad for $l=0$  or $l=1$ for the sequence $ A_n $ the condition  (\ref{f2.1}) is satisfied;

(ii) \quad for $l>1$  the condition (\ref{f2.1}) is equivalent to the condition
\begin{equation} \label{f2.5}
\sum_{k=1}^{\infty} \exp(-R(k+l)+ R(k+1))=\infty.
\end{equation}

\end{lemma}

For the sequence of random events $(A_n)$ we introduce the following
notation
\[
S_{1,n}=  \sum_{k=1}^{n}\mathbf{P}(A_k),  \quad  S_{2,n}=
\sum_{1\leq j < k \leq n}\mathbf{P}(A_j \bigcap A_k).
\]

\begin{lemma}\label{l2.3}
Suppose that in the conditions and notation of Lemma \ref{l2.2}  for some
$l\geq 0$  the condition  (\ref{f2.5}) is satisfied, then
there exists $K>0$ such that
\begin{equation} \label{f2.6}
\limsup_{n\rightarrow\infty} \frac{S_{2,n}}{|S_{1,n}|^2}\leq K <
\infty.
\end{equation}

\end{lemma}

The proofs of Lemmas \ref{l2.2}, \ref{l2.3} can be found in article \cite{mi_tims19}.

 Let us proceed directly to the proof of Theorem \ref{t1}. In the case of
$l=0$ or $l=1$  according to item (i) of Lemma \ref{l2.2}  the series
$\sum_{n=1}^{\infty}\mathbf{P}(A_n )$ diverges, where random
events $A_n$ are given by the equality  (\ref{f2.4}).

Let $l>1 $.  Suppose that the series  (\ref{f7}) diverges for $j=l-1$.
Since under the conditions of Theorem \ref{t1} $r(n)$  is a non-decreasing function, then
\begin{equation} \label{f2.7}
(l-1)r(k+1) \leq R(k+l) - R(k+1)   \leq (l-1)r(k+l).
\end{equation}
 That is, the condition  (\ref{f2.5}) is equivalent to
condition:
\begin{equation}\label{f2.8}
\sum_{k\geq 1}  \exp(- (l-1) r(k)) = \infty .
\end{equation}
 This means that, under the conditions of Theorem  \ref{t1}, in the case where $l>1$, the series
$\sum_{n=1}^{\infty}\mathbf{P}(A_n )$ also diverges.

Then, according to Lemmas \ref{l2.2},\ref{l2.3}, from condition
(\ref{f2.8}) for  $l \geq 0 $ we obtain the estimate (\ref{f2.6}). If
we also take into account the following equality
\begin{equation}\label{f2.80}
\sum_{j=1}^{n} \sum_{l=1}^{n}\mathbf{P}(A_j \bigcap A_l)= 2 S_{2,n}
+S_{1,n},
\end{equation}
then we also have inequalities (\ref{f2.2}), (\ref{f2.3}). It
remains to use Lemma  \ref{l2.1} and Remark \ref{r2}. Thus, for random
 events  $ A_n $, given by the equality  (\ref{f2.4}),
we obtain:
\[
\mathbf{P}(A_n \quad \mbox{occur  \quad i.o.})=1.
\]
Since
\[
A_n \in (\bar{\xi}_n = a_n +l) \quad \mbox{a.s.},
\]
then the following equality immediately follows
 \[ \mathbf{P}(\bar{\xi}_n = a_n +l
\quad \mbox{i.o.} ) = 1.
 \]

On the contrary, suppose that for  $l > 1$
\begin{equation}\label{f2.9}
\sum_{k\geq 1}  \exp(- (l-1) r(k)) < \infty.
\end{equation}
It can be shown (see the proof of Theorem 1  in \cite{mi16}) that
\begin{eqnarray}\label{f2.10}
&  &  \sum_{n\geq 1} \mathbf{P}(\xi \geq   a_n +l )
  =
\sum_{k\geq 1} \mathbf{P}(\xi \geq k +l ) \exp(R(k+1))(1+o(1))
 \nonumber   \\
 & = &
\sum_{k\geq 1} \exp(-R(k+l)+R(k+1))(1+o(1))  .
 \end{eqnarray}
Putting together the estimates (\ref{f2.7}), (\ref{f2.9}),
(\ref{f2.10}) we get
 \[  \sum_{n\geq 1} \mathbf{P}(\xi \geq   a_n +l )< \infty .
\]
As it is well known (see \cite{gal}, $\S 4.3$, Corollary 4.3.1) the last estimate implies the following equality
 \[
 \mathbf{P}(\bar{\xi}_n \geq  a_n +l \quad \mbox{i.o.})= 0 .
\]
So
 \[
 \mathbf{P}(\bar{\xi}_n =  a_n +l \quad \mbox{i.o.})= 0,
\]
which completes the proof of the Theorem  \ref{t1}.
 \quad $ \Box $

Corollary  \ref{c1}  follows directly from Theorem \ref{t1}.  And
equality (\ref{f5}) and  equality (\ref{f8}) for $l=-1$
are simple  consequences of item (ii) of Theorem A and condition
(\ref{f3}),  since for sufficiently large  $k  \quad$ $m \cdot r(k) \leq \exp(r(k))$.

\begin{remark}\label{r3}
If under the condition of Theorem \ref{t1} the series (\ref{f7})
diverges for all $j \geq 0$, then for any $l \geq 0$  the equality
(\ref{f9}) holds.

\end{remark}

\section{Proof of Theorem \ref{t2}.}

The proof of Theorem \ref{t2} will largely be based on the following
important Lemma that was established in \cite{{kla85}}.

\begin{lemma}\label{l3.1}
Let $\xi , \xi_1 , \xi_2 , \ldots ,$ be a sequence of i.i.d.r.v.,
 $(u_{n})$ is  a non-decreasing sequence of real numbers such that, \\ if $n \rightarrow
\infty$
\[
\mathbf{P}( \xi > u_{n}) \rightarrow 0 , \quad n\mathbf{P}( \xi >
u_{n}) \rightarrow \infty,
\]
then the probability
$$\mathbf{P}\left(\bar{\xi}_{n}\leq u_{n} \quad \mbox {i.o.} \right)$$
equals zero or one, depending on whether the following series converges or
diverges
\begin{equation} \label{f3.1}
\sum_{n=1}^{\infty}\mathbf{P}( \xi > u_{n})\exp\left(-n\mathbf{P}(
\xi > u_{n})\right) .
\end{equation}
Furthemore, if
$$\mathbf{P}( \xi > u_{n}) \rightarrow c > 0,$$
then
$$\mathbf{P}\left(z_{n}\leq u_{n} \quad \mbox {i.o.}
\right) =0 ,$$
and if
$$\liminf_{n \rightarrow \infty} n\mathbf{P}( \xi >
u_{n}) <\infty,$$
then
$$\mathbf{P}\left(z_{n}\leq u_{n} \quad
\mbox {i.o.} \right) =1 .$$

\end{lemma}

We will use the following notations from now on:
\[
I_k = \{n\geq 1 : a_{n} =k \} = \{n\geq 1 : \exp(R(k)) \leq n <
\exp(R(k+1)) \} , \quad k = 1, 2, \ldots ,
\]
\[
d_k =\mathbf{P}(\xi > k - l )  = \exp(-R(k-l+1)), \quad s_k =
\sum_{n \in I_k }  \exp(-n d_k ).
 \]

(i) \quad Let  $l=1$,  $u_n =a_n -1$,  $n\in I_k$.  We choose  $n$ such that
 $n=\exp(R(k))+ \delta$, $0 \leq \delta \leq 1$ .

Then, for sufficiently large $k$
 \[
 n\mathbf{P}( \xi > u_{n})= n\mathbf{P}( \xi \geq k)
 = (\exp(R(k))+ \delta) \exp(-R(k)) \sim 1 +o(1),
 \]
 that is
\[
\liminf_{n \rightarrow \infty} n\mathbf{P}( \xi > u_{n}) < \infty .
 \]
By Lemma \ref{l3.1} this implies the equality  (\ref{f10}) for
$l=1$.

(ii) \quad Let $l>1 $ \, and
\begin{equation} \label{f3.111}
\sum_{k\geq 1}  \exp\left(- \exp((l-1) r(k))\right) < \infty .
\end{equation}
Let us estimate the series (\ref{f3.1}) from Lemma \ref{l3.1} for $u_n = a_{n} - l $:

\begin{eqnarray}\label{f3.2}
&  &  \sum_{n\geq 1} \mathbf{P}(\xi >   a_n - l )  \exp(-n
\mathbf{P}(\xi >   a_n - l) ) \nonumber   \\
 & = &
\sum_{k\geq 1} \mathbf{P}(\xi > k - l ) \sum_{n \in I_k } \exp(-n
\mathbf{P}(\xi >   k - l ) )  \nonumber   \\
 & = &
\sum_{k\geq 1} d_k  s_k  .
 \end{eqnarray}

To estimate the last series, we use the following known
inequalities:\\ if $f(t)$ is a nonincreasing function, then for  $i\leq j$
\begin{equation} \label{f3.3}
\sum_{k =i+ 1}^{j+1} f(k) \leq \int_{i}^{j+1} f(t) dt \leq \sum_{k
=i}^{j} f(k) .
\end{equation}
Let's substitute into the left-hand side of inequality  (\ref{f3.3})
 function $ f(t)=\exp(-td_k )$. This allows us to write an upper bound for $s_k$.

\begin{eqnarray} \label{f3.4}
s_k & \leq & \sum_{[\exp(R(k)]\leq n \leq [\exp(R(k+1)]}  \exp(-n
d_k) \leq \int_{\exp(R(k))-2}^{\exp(R(k+1))} \exp(-t d_k ) dt \nonumber \\
& \leq &
 \frac{1}{d_k }\left( \exp(-d_k  (e^{R(k)} -2)) - \exp(-d_k e^{R(k+1)} )\right).
\end{eqnarray}

Consider separately the first term in brackets from the
estimate (\ref{f3.4}):
\begin{eqnarray} \label{f3.5}
 \exp\left(-d_k  (e^{R(k)} -2)\right) & = & \exp\left(-e^{R(k- l+1)} (e^{R(k)} -2)\right) \nonumber \\
= \exp\left(-e^{r(k)+ \ldots +r(k-l+2)} +o(1)\right) & \leq & \exp\left(-e^{(l-1)
r(k-l+2)})(1 +o(1)\right).
\end{eqnarray}
 The second term  is equally easy to estimate
\begin{eqnarray} \label{f3.6}
 \exp\left(-d_k  e^{R(k+1)}\right) & = & \exp\left(-e^{-R(k- l+1)} e^{R(k+1)} \right) \nonumber \\
= \exp\left(-e^{r(k+1)+ \ldots +r(k-l+2)} \right) & \leq & \exp\left(-e^{l
r(k-l+2)}\right) .
\end{eqnarray}
 Together, the
estimates   (\ref{f3.4})-(\ref{f3.6}) give
\begin{eqnarray*}
 \sum_{k\geq 1} d_k  s_k \leq \sum_{k\geq 1} \exp(-e^{(l-1)
r(k-l+2)})(1+ o(1)).
\end{eqnarray*}
 The convergence of the last series is equivalent to the
condition (\ref{f3.111}). Considering the equality
 (\ref{f3.2}), we also obtain the convergence of series (\ref{f3.1}).

The other conditions of Lemma  \ref{l3.1}  also hold. For example,
let $u_n =a_n -l$, $l>1$, $n\in I_k$,  that is $\exp(R(k)) \leq n < \exp(R(k+1))$. Then
 \[
 n\mathbf{P}( \xi > u_{n})= n\mathbf{P}( \xi > k-l)
 \geq  \exp(R(k)) \exp(-R(k-l+1))
 \]
\[
 \geq  \exp((l-1)r(k-l))
 \rightarrow \infty , \quad k \rightarrow \infty .
 \]
Thus, the probability (\ref{f111})  is equal to zero.

 In the next step, let's assume that
\begin{equation} \label{f3.7}
\sum_{k\geq 1}\exp\left(-e^{(l-1)r(k)}\right) = \infty,
\end{equation}
and show that
\begin{equation} \label{f3.8}
\mathbf{P}(\bar{\xi}_n \leq a_n -l  \quad \mbox{i.o.} )=1.
\end{equation}

To do this, we need to find a lower bound for the series  (\ref{f3.1})
from Lemma \ref{l3.1}  for $u_n = a_{n} - l  $.  It is clear
that the equality (\ref{f3.2}) will remain unchanged. To find a lower estimate for the value ${s}_k $, we again use the inequalities (\ref{f3.3}),  but now we choose the right-hand inequality.
\begin{eqnarray} \label{f3.9}
{s}_k & \geq & \sum_{[\exp(R(k)]+1 \leq n \leq [\exp(R(k+1)]-1}
\exp(-n
{d}_k) \geq \int_{\exp(R(k))+1}^{\exp(R(k+1))-1} \exp(-t {d}_k ) dt \nonumber \\
& \leq &
 \frac{1}{{d}_k }\left( \exp(-{d}_k  (e^{R(k)} +1)) - \exp(-{d}_k e^{R(k+1)-1} )\right).
\end{eqnarray}
It should be noted that
\begin{eqnarray} \label{f3.10}
  \exp(-{d}_k  \left(e^{R(k)} +1)\right) &=& \exp\left(- e^{r(k)+ \ldots +r(k-l+2)}+o(1)\right) \nonumber \\
   & \geq &
     \exp(- e^{(l-1)r(k)})(1+o(1))
\end{eqnarray}
and
\begin{eqnarray} \label{f3.11}
  \exp\left(-{d}_k e^{R(k+1)-1} \right) &=& \exp\left(- e^{r(k+1)+ \ldots +r(k-l+2)}+o(1)\right) \nonumber \\
  &=& o(
     \exp\left(- e^{r(k)+ \ldots +r(k-l+2)})\right).
\end{eqnarray}

Thus, according to the relations (\ref{f3.9}) - (\ref{f3.11}) and
taking into account the condition (\ref{f3.7}),  we get
\[
\sum_{k\geq 1} {d}_k  {s}_k \geq  \sum_{k\geq 1} \exp\left(-
e^{(l-1)r(k)}\right)(1+o(1))=\infty .
\]
By Lemma  \ref{l3.1} this means that the equality(\ref{f3.8}) is
true. $\Box $

\section{Proof of Theorem \ref{t3}.}

To prove the equality (\ref{f100})  under the condition
\begin{equation} \label{f4.1}
\sum_{k\geq 1}\exp\left(-e^{(l-1) \, r(k)}\right) = \infty
\end{equation}
we will use some generalization of the method of the
article  \cite{adm24}.

For this purpose let's introduce random events  $\hat{A}_k$, $k\geq l$:
\begin{equation}\label{f4.2}
\hat{A}_k = \hat{A}_{k,l} = \{\bar{\xi}_{n_k} =k-l     \},
\end{equation}
where
\[
n_k =\min(n\geq 1:  a_n = k) = \lceil  \exp(R(k)) \rceil , \quad k
\geq l,
\]
\quad where $\lceil x \rceil$  - is the smallest integer $\geq x$.

Since
\[
 \quad
\{\bar{\xi}_{n_k} =k-l \quad \mbox{i.o.} \}  \subset \{\bar{\xi}_{n}
=a_n-l \quad \mbox{i.o.}\},
\]
  then equality (\ref{f100}) is
a simple consequence of the following relation
\begin{equation}\label{f4.3}
 \quad   \mathbf{P}\left(\bigcap_{n=1}^{\infty} \bigcup_{k=n}^{\infty} \hat{A}_k \right)
 =1.
\end{equation}

As in the case of Theorem \ref{t1},  we use Lemma \ref{l2.1} to
prove the equality (\ref{f4.3}).

Step 1. Let us establish the following equality:
\begin{equation}\label{f4.4}
 \sum_{k\geq 1}  \mathbf{P}(\hat{A}_k )= \infty .
\end{equation}

Directly from the definition we get
 \begin{eqnarray}\label{f4.17}
 \mathbf{P}(\hat{A}_k) & = &
 \mathbf{P}( \bar{\xi}_{n_k} \leq k-l ) - \mathbf{P}( \bar{\xi}_{n_k} \leq k-l-1 )
  \nonumber \\
 & = &
(1- \mathbf{P}( {\xi}\geq k-l+1 ))^{n_k} - (1- \mathbf{P}(
{\xi}\geq k-l))^{n_k}
 \nonumber \\
 & = &
(1- \exp(-R(k-l+1) ))^{n_k} - (1- \exp(-R(k-l) ))^{n_k}
\nonumber \\
 & = &
V_1 (k) \cdot (1- V_2 (k)),
  \end{eqnarray}
where
\begin{eqnarray}\label{f4.18}
  V_1(k) &=& \left((1- \exp(-R(k-l)+1))\right)^{n_k},  \nonumber \\
 V_2(k) &=&   \left(\frac{1- \exp(-R(k-l))}{1- \exp(-R(k-l+1))}\right)^{n_k}.
  \end{eqnarray}

It is well known that  $(1-1/x)^x \uparrow 1/e$  for $x
\uparrow \infty$. And since $(1-1/x)^{x-1} $ is a decreasing
function for $x > 1 $,  then

\begin{eqnarray}\label{f4.19}
\left(1-\frac{1}{x}\right)^{x-1}  \downarrow \frac{1}{e} \quad \mbox{for} \quad  x
\uparrow \infty.
 \end{eqnarray}

Then, based on the relation (\ref{f4.19})  for sufficiently large $k$ we get
\begin{eqnarray}\label{f4.20}
  V_1(k) &=&  (1- \exp(-R(k-l+1)))^{ \lceil
\exp(R(k)) \rceil }
\nonumber \\
 & = &
(1- \exp(-R(k-l+1)))^{  \exp(R(k))+\theta}
   =  V_{1,1} (k)\cdot V_{1,2} (k),
  \end{eqnarray}
where \quad $0 \leq \theta \leq 1$,
\begin{eqnarray}\label{f4.21}
  V_{1,1} (k) &=&  \left((1- \exp(-R(k-l+1)))^{
\exp(R(k-l+1)) -1}\right)^{\exp(R(k)-R(k-l+1))}
\nonumber \\
 & \geq &
\left(\frac{1}{e}\right)^{\exp(R(k)-R(k-l+1))}
   \geq \exp(-\exp((l-1) \cdot r(k)));
  \end{eqnarray}
\begin{eqnarray}\label{f4.22}
  V_{1,2} (k) &=&  \left(1- \exp(-R(k-l+1))\right)^{
\exp(R(k)-R(k-l+1))+\theta}
\nonumber \\
 & \geq & (1-\epsilon)
\left((1- \exp(-R(k-l+1)))^{ \exp(R(k-l+1))}\right) ^{\exp(R(k)-2 R(k-l+1))}
  \end{eqnarray}
for a small $\epsilon >0$.

Under the conditions of Theorem  \ref{t3} the function $r(k)$
increases and $r(k) = o(k)$, and  $R(k)\geq r(1) k$.  Then, for $ k\rightarrow\infty $
\[
R(k) - 2R(k-l+1)\leq (l-1)r(k) -r(1) (k-l+1 )
 \rightarrow -\infty,
\]
\[
\left(1- \exp(-R(k-l+1))\right)^{ \exp(R(k-l+1))} \rightarrow  \frac{1}{e} .
\]
 Thus $ V_{1,2} (k) \geq 1- \epsilon$ for $k\rightarrow\infty$.

 From this and the relations  (\ref{f4.20})-(\ref{f4.22}) we get:
 for any  $\epsilon > 0$  with sufficiently large $ k$ we get:
\begin{eqnarray}\label{f4.23}
 V_1 (k)  \geq (1-\epsilon)
 \exp(-\exp((l-1) \cdot r(k))).
  \end{eqnarray}

Next, we find the upper bound for the value $V_2 (k)$:
\begin{eqnarray}\label{f4.24}
V_2 (k)  & = &  \left(1- \frac{ \exp(-R(k-l)) - \exp(-R(k-l+1))}{1-
\exp(-R(k-l+1))}\right)^{n_k}
\nonumber \\
 & = &
  \left(1- \frac{ \exp(-R(k-l))(1- \exp(-r(k-l+1)))}{1-
\exp(-R(k-l+1))}\right)^{n_k}
\nonumber \\
 & = &
 \left(\left(1-\frac{1}{y}\right)^y \right)^{\frac{n_k}{y}},
  \end{eqnarray}
where
\[
y= \frac{\exp(R(k-l))(1- \exp(-R(k-l+1))) }{ 1- \exp(-r(k-l+1))} .
\]

Since
\[
\frac{n_k }{y}= \frac{\lceil \exp(R(k)) \rceil }{y} \rightarrow
\infty,
\]
then
\begin{eqnarray}\label{f4.25}
V_2 (k)= \left(\left(1-\frac{1}{y}\right)^y \right)^{\frac{n_k}{y}}\leq
\left(\frac{1}{e}\right)^{\frac{n_k}{y}} \rightarrow 0.
\end{eqnarray}
Taking together the estimates (\ref{f4.23}), (\ref{f4.25}) and
(\ref{f4.17})  considering the condition (\ref{f4.1}), we obtain the
equality  (\ref{f4.4}).

Step 2. To apply Lemma \ref{l2.1} to random events  $\hat{A}_n$
 we need to prove the following inequality:
\begin{equation} \label{f4.230}
\limsup_{n\rightarrow\infty} \frac{\hat{S}_{2,n}}{|\hat{S}_{1,n}|^2}
\leq K < \infty,
\end{equation}
where
\[
\hat{S}_{1,n}=  \sum_{k=1}^{n}\mathbf{P}(\hat{A}_k),  \quad
\hat{S}_{2,n}= \sum_{1\leq j < k \leq n}\mathbf{P}\left(\hat{A}_j \bigcap
\hat{A}_k\right).
\]
Indeed, this will be clear if we consider equalities  (\ref{f4.4}) and  (\ref{f2.80}).

 Let us write the value  $\hat{S}_{2,n}$ in the following form
\begin{equation} \label{f4.240}
 \hat{S}_{2,n}= \sum_{1\leq j < k \leq n}\mathbf{P}(\hat{A_j}) \mathbf{P}(\hat{A}_k / \hat{A}_j).
\end{equation}
Next we find the upper bound of another factor in (\ref{f4.240})
\begin{eqnarray}\label{f4.26}
  \mathbf{P}(\hat{A}_k / \hat{A}_j)& =&
  \mathbf{P}\left(\bar{\xi}_{n_k} =k -l
\ldots /\bar{\xi}_{n_j} =j -l\right)
 =   \mathbf{P}\left(\max_{s\in (n_j +1, n_k)} \xi_s =k -l\right)
\nonumber \\
& = &  \mathbf{P}\left( \max_{s\in (n_j +1, n_k)} \xi_s =k -l,  \quad
\xi_i \leq k -l, \quad i=1, \ldots, n_j, \right) C_{j,k}
 \nonumber \\
& \leq &
 \mathbf{P}(\hat{A}_k )  C_{j,k},
   \end{eqnarray}
where  \quad $C_{j,k}= (\mathbf{P}( \xi_i \leq k -l, i=1, \ldots,
n_j ))^{-1}$.

Since for $x>1$ the function  $(1 - {1}/{x}
)^{-x}$ decreases as  $x$ increases, then for  $ 1\leq j
\leq k-l $
\begin{eqnarray}\label{f4.27}
C_{j,k} & = & (1-\exp(-R(k-l+1)))^{-\lceil \exp(R(j)) \rceil}
 \nonumber \\
& \leq & (1-\exp(-R(k-l+1)))^{- \exp(R(k-l+1)) -1 }
\nonumber \\
& \leq & (1-\exp(-R(1)))^{- \exp(R(1)) -1 } = C_0 .
   \end{eqnarray}

Therefore
\begin{equation} \label{f4.28}
\sum_{1 \leq k \leq n} \sum_{1 \leq j\leq k-l }\mathbf{P}\left( \hat{A}_j
\bigcap \hat{A}_k\right) \leq  C_0 \sum_{1 \leq k \leq n} \sum_{1 \leq
j\leq k-l }\mathbf{P}( \hat{A}_j ) \mathbf{P}( \hat{A}_k) .
\end{equation}

The case for  $k-l < j \leq k $ is also investigated simply
\begin{equation} \label{f4.29}
\sum_{1 \leq k \leq n} \sum_{k-l  < j \leq k }\mathbf{P} \left( \hat{A}_j
\bigcap \hat{A}_k\right) = \sum_{1 \leq k \leq n} \sum_{k-l  < j \leq k
}\mathbf{P}( \hat{A}_k ) \mathbf{P}( \hat{A}_j /\hat{A}_k) \leq l
\sum_{1 \leq k \leq n}\mathbf{P}( \hat{A}_k).
\end{equation}

Finally, from equality (\ref{f4.4}) we get: $\hat{S}_{1,n}
\rightarrow \infty $ for  $n \rightarrow \infty $
 and
\begin{equation} \label{f4.30}
 (\hat{S}_{1,n})^2 = 2 \sum_{1 \leq j < k \leq n} \mathbf{P}( \hat{A}_j ) \mathbf{P}(
 \hat{A}_k ) + O(\hat{S}_{1,n}) .
\end{equation}

Putting together the estimates  (\ref{f4.28}) - (\ref{f4.30}) we
obtain

\[
\hat{S}_{2,n}\leq  C_0 \sum_{1 \leq k \leq n} \sum_{1 \leq j\leq k-l
}\mathbf{P}( \hat{A}_j ) \mathbf{P}( \hat{A}_k) + l \sum_{1 \leq k
\leq n}\mathbf{P}( \hat{A}_k)
\]
\[
\leq  C_0 (\hat{S}_{1,n})^2  + l  \hat{S}_{1,n},
\]
that gives the estimate (\ref{f4.230}).

Hence, from Lemma \ref{l2.1}  and Remark $2$ we immediately
obtain equality (\ref{f100}). $\Box $

Corollary \ref{c2} simply follows from Theorems \ref{t2}, \ref{t3}.
Indeed, by Theorem  \ref{t2} we get
\[
\mathbf{P}(\bar{\xi}_n \leq  a_n -m -1 \quad \mbox{i.o.} ) = 0,
\quad \mathbf{P}(\bar{\xi}_n \leq  a_n -m  \quad \mbox{i.o.} ) = 1 .
\]
 Hence, we obtain the equality  (\ref{f13}).

Item (ii)  is a direct consequence of Theorems   \ref{t1} and \ref{t3} .

\section{ Examples}

Let us give some examples of using the obtained results to study the asymptotics of extreme values.

\emph{Example 1. Poisson distribution}. Let $\xi$ be a r.v. with the distribution \quad $p_i = \frac{\lambda^i}{i!}
\exp(-\lambda), i\geq 0$.  As it was established in \cite{mi16} for
the Poisson distribution the following asymptotic formulas are true
\[
R(n)= \left(n+\frac{1}{2}\right)\ln n -n(\ln \lambda +1)+ \frac{1}{2} \ln 2\pi
- \lambda +o(1),
\]
\[
r(n)= \ln n + O(1).
\]
It is easy to verify that the conditions of Corollary \ref{c1}  hold
for $m=2$. Therefore, for it the following equalities are true:
  \[
  \mathbf{P}\left(\limsup_{n \rightarrow\infty} (\bar{\xi}_n -  a_n ) = 2 \right)
= 1,
\]
 \[
  \mathbf{P}\left(\liminf_{n \rightarrow\infty} (\bar{\xi}_n -  a_n ) = -1 \right)
= 1,
\]
and also the quality  (\ref{f9}) for $ l=-1, 0, 1,2 $, where
\begin{equation*}
a_n =   \frac{\ln n}{ L_2 (n)}\left(1+ \frac{L_3 (n)+\ln \lambda  + 1
+o(1)}{L_2 (n)}\right)
\end{equation*}
(see also  \cite{mi_tims19}).

\emph{Example 2.}  Suppose $\xi$ is a discrete r.v., for which the
values  $r(n)$ and $a(n)$  are given by equalities  (\ref{f2}) and
for some  $C$, $0 < C < 1$
\[
r(n) =(C +o(1)) \cdot L_2 (n)\uparrow \quad \mbox{for} \quad n
\uparrow \infty .
\]
Let's put $m=\min(l\geq 1 : l \cdot C>1)$.  Then the series
(\ref{f12}) in Corollary  \ref{c2} \, converges for $j=m $,
 and diverges for $j=m-1 $.  Therefore, for extreme values of
independent copies of a given  r.v. $\xi$ the equality (\ref{f13})
is satisfied, as well as the equality (\ref{f100})   $\forall \, l \in (m, m-1, \ldots,
1)  $.

\emph{Example 3. Geometric distribution.} Let  $\xi$ be a r.v. with
the distribution $ \mathbf{P}(\xi =i) = p_i =q(1-q)^i , i\geq
0$ . Then
\[
\mathbf{P}(\xi \geq n )=(1-q)^n =\exp(-\gamma n), \quad \gamma =\ln
1/(1-q).
\]
It is easy to find
\begin{equation} \label{f146}
R(n)  =\gamma n, \quad r(n) =\gamma,  \quad a_n = [(1/\gamma) L_{1}
(n) ],
\end{equation}
wher $[x]$ is an integer part of the number $x$.
It is clear that the geometric distribution satisfies the condition
of item (ii) of Theorem B.  Therefore, for any integer $m$, the
equality  (\ref{f130}) is true.

\begin{remark}\label{r4}
 A random variable with a geometric distribution can be
represented as follows.  Let $ \tau^e $ be an exponential r.v. with
the parameter   $\gamma > 0 $, that is  $ \mathbf{P}(\tau^e < x
)=1- \exp(-\gamma x)$. Let us consider a r.v. $ \xi = [\tau^e]
$.  Then
\[
\mathbf{P}(\xi = k )=\mathbf{P}(\tau^e \in [k, k+1] )= q (1-q)^k ,
\quad k=0, 1, 2, \ldots ,
\]
where $q= 1- \exp(-\gamma )$, that is  $ \xi $ has a
geometric distribution. And if  $\bar{\tau}^e_n $ and
$\bar{\xi}_n$  are corresponding sequences of maxima,
constructed from r.v.  $ \tau^e $ and ${\xi}$, then
\[
 \bar{\tau}^e_n = \bar{\xi}_n + \theta_n ,  \quad 0 \leq \theta_n \leq
 1  \quad a.s.
\]
Hence and Theorem B for any fixed integer $m$ we get
\[
\mathbf{P}(\bar{\xi}_n = a_n +m \quad \mbox{i.o.})=
\mathbf{P}(\bar{\tau}^e_n \in [a_n +m, a_n +m +1]   \quad
\mbox{i.o.} )= 1 .
\]
where $a_n$ is defined in equality (\ref{f146}) .

It is clear that Theorems \ref{t1}-\ref{t3} can also be extended
to arbitrary continuous r.v.

 \end{remark}

\vspace{0.5cm}


\end{document}